\documentclass[12pt]{article}

\usepackage{listings}

\usepackage{epsfig}
\usepackage{epsf}
\usepackage{amsmath}
\usepackage{amsfonts}
\usepackage{amsthm}
\usepackage{amssymb}
\usepackage{enumerate}

\usepackage{pstricks}
\usepackage{lastpage}

\def\C{\mathbb{C}}
\def\R{\mathbb{R}}
\def\i{\bf i}

\usepackage{xspace}
\usepackage{tikz} 
\usepackage{calc} 
\usepackage{ifthen} 

\usetikzlibrary{decorations.pathmorphing, decorations.pathreplacing, patterns, shapes.geometric} 

\tikzset{
    slope/.code={\edef\slope{#1}},
    slope/.default=0.5,
    slope
}
\makeatletter
\pgfdeclarepatternformonly[\tikz@pattern@color,\slope]{slant lines}
{\pgfpoint{-.1mm/\slope}{-.1mm}}
{\pgfpoint{1.1mm/\slope}{1.1mm}}
{\pgfpoint{1mm/\slope}{1mm}}
{
    \pgfsetlinewidth{0.4pt}
    \pgfpathmoveto{\pgfpoint{-.1mm/\slope}{-.1mm}}
    \pgfpathlineto{\pgfpoint{1.1mm/\slope}{1.1mm}}
    \pgfsetstrokecolor{\tikz@pattern@color}
    \pgfusepath{stroke}
}
\makeatother

\usepackage{fancyhdr} 
\begin{document}

\baselineskip=18pt

\begin{center}
{\large\bf

THE MULTI-DIMENSIONAL DECOMPOSITION WITH CONSTRAINTS

}

\smallskip

Ilgis Ibragimov, Elena Ibragimova

\smallskip

{\it

Elegant Mathematics LLC, 82834 WY USA \&

Elegant Mathematics Ltd, 66564 Germany

\smallskip

e-mail: ii@elegant-mathematics.com

}

\end{center}

\bigskip

{\bf\large

\noindent
ABSTRACT

}
\bigskip

We search for the best fit in the Frobenius norm of $A \in \C^{m \times n}$
by a matrix product $B C^*$, where $B \in \C^{m \times r}$ and
$C \in \C^{n \times r}$, with  $r \le m$ so that $B = \{b_{ij}\}_{\tiny
\begin{tabular}{l}i=1, \dots, m \\ j=1, \dots, r \end{tabular}}$
is defined by some unknown parameters $\sigma_1, \dots, \sigma_k$,
$k << mr$, and all partial derivatives of $\displaystyle
\frac{\delta b_{ij}}{\delta \sigma_l}$ are definite, bounded, and can
be computed analytically.

We show that this problem transforms to a new minimization problem
with only $k$ unknowns by the analytical computation of the gradient of
the minimized function over all $\sigma$. The complexity of computation of
this gradient is only 4 times greater than the complexity of computation of
the function, and this new algorithm needs only $3mr$ additional words in memory.

We apply this approach for the solution of the three-way decomposition problem
and obtain good result of convergence for the Broyden algorithm.

\bigskip

{\bf\large

\noindent
INTRODUCTION

}
\bigskip

Suppose we have $A \in \C^{m \times n}$. The idea is to find
$B \in \C^{m \times r}$ and $C \in \C^{n \times r}$, $r \le m$ so

\begin{equation}
\label{eq:1}
\min_{C, \sigma_1, \dots, \sigma_k} || A - B(\bar \sigma) C^*||_F^2,
\end{equation}

\noindent that $B = \{b_{ij}\}_{\tiny \begin{tabular}{l}i=1, \dots,
m \\ j=1, \dots, r \end{tabular}}$ is defined by unknown parameters
$\sigma_1, \dots, \sigma_k$, $k << mr$, and all partial derivatives
of $\displaystyle \frac{\delta b_{ij}}{\delta \sigma_l}$ are definite,
bounded and can be computed analytically.

This problem occurs in statistics \cite{St}, nuclear magnetic resonance
\cite{NMR}, spectroscopy and multi-dimensional decomposition \cite{ENUMATH}.
Consider one popular application \cite{arXiv-2} --- a low rank approximation
of two and multidimansional data array with one factor matrix containing
vectors formed as complex exponents:
\begin{equation}
\label{eq:eqq1}
\min_{B, \sigma} \sum_{j=1}^J \sum_{k=1}^K \left|\left| a_{jk} - \sum_{l=1}^L b_{jl} e^{\i \sigma_{k l}} \right|\right|_2^2, 
\end{equation}
and
\begin{equation}
\label{eq:eqq2}
\min_{B, \sigma} \sum_{j=1}^J \sum_{k=1}^K \left|\left| a_{jk} - \left| \sum_{l=1}^L b_{jl} e^{\i \sigma_{k l}} \right| \right|\right|_2^2, 
\end{equation}
%
Since the total amount of minimizing paramenters $\sigma$ usually is several orders less than the
total amount of minimizing paramenters in $B$, it is highly desired to perform
minimization over only $\sigma$ to save computational complexity.

If we freeze $B$, then this function is linear in $C$, and
$C = A^* B (B^* B)^{-1}$. The problem (\ref{eq:1}) then turns into a new
nonlinear problem with only $k$ unknowns:

\begin{equation}
\label{eq:2}
\min_{\sigma_1, \dots, \sigma_k} || A - B (B^* B)^{-1} B^* A||_F =
\min_{\sigma_1, \dots, \sigma_k} \sqrt{|| A ||_F^2 - ||A^* Q(B)||_F^2},
\end{equation}

\noindent where $Q(B) \in \C^{m \times r}$ contains the orthonormal subspace
from $B$.

The main difficulty in applying minimization methods for (\ref{eq:2})
is the computation of the gradient of the function over all $\sigma$. The finite
difference method needs $k$ or $2k$ computations of this function
for one evaluation of the gradient and cannot be considered accurate.
There is a good alternative for it, Baur-Strassen (BS) method
\cite{BS}, which allows computing the gradient of a function using only $5n$
operations if the original function can be computed by $n$ simple
arithmetical operations with no more than 2 operands. The big disadvantage
of the BS method is its memory requirement: it needs ${\cal O}(n)$ words
in memory, which is too many for most applications.

We suggest a new approach for computing the gradient of a function.
This approach contains Modified Gramm--Schmidt (MGS) orthogonalization
with low memory requirements and is based on the BS method.

\bigskip

{\bf\large

\noindent
ALGORITHM

}
\bigskip

To compute (\ref{eq:2}), we perform the following steps:

\begin{itemize}

\item [{\bf 1)}] create $B$ from $\sigma_1, \dots, \sigma_k$;

\item [{\bf 2)}] compute orthonormal subspace $Q$ in $B$;

\item [{\bf 3)}] compute (\ref{eq:2}).

\end{itemize}

\noindent In this article, we discuss how to compute a gradient
$\hat g \in \C^{mr}$ of (\ref{eq:2}) over all entries of $B$. We
will use both $G \in \C^{m \times r}$ and $\hat g$ for the same data.
Let the dependence of $B$ on $\sigma_1, \dots, \sigma_k$ be so simple
that one can compute the gradient of (\ref{eq:2}) by $\sigma_1, \dots,
\sigma_k$ if $G$ is known.

Steps 2 and 3 need $mr$ additional words in memory and compute within
$2mr(r+n)$ arithmetical operations in the event that the MGS algorithm is used
for step $2$. The BS algorithm can compute the gradient with the same
order of arithmetical complexity but needs $4mr(r+n)$ additional words
in memory.

Let us consider a computation of (\ref{eq:2}) from $B$. Let $B=[b_1,
\dots, b_k]$ be the initial matrix and $Q=[q_1, \dots, q_k]$ the
orthonormal subspace, which we are going to compute. Then

$\displaystyle q_1 = \frac{b_1}{||b_1||_2}$ \\
$do$ $i=2$, $r$ \\
\hspace*{10mm} $u=b_i,$ \\
\hspace*{10mm} $do$ $j=1$, $i-1$ \\
\hspace*{20mm} $u = u - q_j q_j^* u$ \\
\hspace*{10mm} $enddo$ \\
\hspace*{10mm} $\displaystyle q_i = \frac{u}{||u||_2}$ \\
$enddo$ \\
$\displaystyle f=\sqrt{|| A ||_F^2 - \sum_{i=1}^r ||A^* q_i||_2^2}$ \\

Let's construct a gradient of $f$ by $B$. We will call ${\bf d} y_i \in \C^m$
the vector of derivatives --- each $k$-th element of this vector contains
the derivative of the $k$-th element of vector $y_i$. Then there are the
following formulas for the gradient:

$\displaystyle {\bf d} q_1 = \frac{1}{||b_1||_2} (I - q_1 q_1^*)
{\bf d} b_1$ \\
$do$ $i=2$, $r$ \\
\hspace*{10mm} $u=b_i$ \\
\hspace*{10mm} $do$ $j=1$, $i-1$ \\
\hspace*{20mm} ${\bf d} u_{new} = (I - q_j q_j^*) {\bf d} u_{old} -
(q_j^* u_{old} I + u_{old} q_j^*) {\bf d} q_j$ \\
\hspace*{10mm} $enddo$ \\
\hspace*{10mm} $\displaystyle {\bf d} q_i = \frac{1}{||u||_2}
                (I - q_i q_i^*) {\bf d} u$ \\
$enddo$ \\
$\displaystyle {\bf d} f = - \frac{1}{f} \sum_{i=1}^r q_i^* A A^* {\bf d} q_i$ \\

\noindent We can write all these equations in matrix notation:

$$\left(\begin{tabular}{ccc}$I_{mr \times mr}$ & 0 & 0 \\
$F_{mr \times \frac{mr(r+1)}{2}}$ &
$L_{\frac{mr(r+1)}{2} \times \frac{mr(r+1)}{2}}$ & 0 \\
0 & $h_{\frac{mr(r+1)}{2}}^*$ & 1 \\ \end{tabular} \right)
\left( \begin{tabular}{c} $\frac{\delta}{\delta b_{ij}}$ \\
$\hat g^*$ \end{tabular} \right) =
\left( \begin{tabular}{c} $I_{mr \times mr}$ \\ 0 \\ \end{tabular} \right)
\hspace*{10mm} {\rm or}$$

\begin{equation}
\label{eq:3}
\left(\begin{tabular}{ccc}$I_{mr \times mr}$ &
$F_{mr \times \frac{mr(r+1)}{2}}^*$ & 0 \\
0 & $L_{\frac{mr(r+1)}{2} \times \frac{mr(r+1)}{2}}^*$ &
$h_{\frac{mr(r+1)}{2}}$ \\
0 & 0 & 1 \\ \end{tabular} \right)
\left( \begin{tabular}{c} $\hat g$ \\ $*$ \\ $\vdots$ \\ $*$ \\
\end{tabular} \right) =
\left( \begin{tabular}{c} $1$ \\ $0$ \\ $\vdots$ \\ $0$ \\
\end{tabular} \right),
\end{equation}

\noindent where $L_{\frac{mr(r+1)}{2} \times \frac{mr(r+1)}{2}}$ is a
lower block triangular matrix with block size $m \times m$. This
matrix has $I_{m \times m}$ blocks on the diagonal. The block matrix
$[F,L]$ contains no more than 3 nonzero blocks in each row (see Fig. 1
for one example with $4$ vectors). We marked with $*$ the elements
that we are not interested in; $\hat g$ is the vector of the gradient of
the original function over all $b_{ij}$. To compute $\hat g$, we solve
the linear system (\ref{eq:3}). If we create $L$ and $F$ matrices,
then we need at least $mr(r+3)$ words to solve $L$.

$$\begin{tabular}{|c|c|c|c|c|c|c|c|c|c|c|c|c|c|c|} \hline
       &$b_1$&$b_2$&$b_3$&$b_4$&$q_1$&$u_{12}$&$q_2$&$u_{13}$&$u_{23}$&$q_3$&$u_{14}$&$u_{24}$&$u_{34}$&$q_4$\\\hline
$q_1   $&$S_1$&  &   &   &$I$&      &   &      &      &   &      &      &      &   \\\hline
$u_{12} $&  &$W_2$&   & &$V_2$& $I$  &   &      &      &   &      &      &      &   \\\hline
$q_2    $&   &   &   &   &   &$S_2$ &$I$&      &      &   &      &      &      &   \\\hline
$u_{13} $&   &  &$W_3$& &$V_3$&     &   & $I$  &      &   &      &      &      &   \\\hline
$u_{23} $&   &   &   &   &   &    &$V_3$&$W_3$ & $I$  &   &      &      &      &   \\\hline
$q_3    $&   &   &   &   &   &      &   &      &$S_3$ &$I$&      &      &      &   \\\hline
$u_{14} $&   &   & &$W_4$&$V_4$&    &   &      &      &   & $I$  &      &      &   \\\hline
$u_{24} $&   &   &   &   &   &     &$V_4$&     &      &   &$W_4$ & $I$  &      &   \\\hline
$u_{34} $&   &   &   &   &   &      &   &      &     &$V_4$&     &$W_4$ & $I$  &   \\\hline
$q_4    $&   &   &   &   &   &      &   &      &      &   &      &      &$S_4$ &$I$\\\hline
\end{tabular}$$

\underline{Figure 1} Shows the matrix $[F, L]$ when $B$ has $4$ vectors,
here $S = \frac{I - uu^*}{||u||_2}$, $V = q^* u I + u q^*$,
$W = qq^* - I$.

We suggest an improvement where we need only $4mr$ words to store some
parts of $L$ and $F$ but still compute the solution. Let's remark that
in the loop for the variable $j$, we update $(i-1)$ times vector $u$.
If we store matrices $B$ and $Q$, we can recompute all updates of $u$
from this loop for particular $i$ with $2(i-1)M$ additional arithmetical
operations and store them in one additional array $T \in \C^{m \times r}$.
Then, during backward substitution we recompute all updates of $u$ only
when we need it. Obviously, this occurs $r-1$ times for all $i=r, \dots, 2$.
All multiplications to matrices $S$, $V$, and $W$ need ${\cal O}(m)$
arithmetical operations. Thus, we need only $B$, $Q$, $T$, and $G$ arrays
with size $m \times r$ for this computation. Here is an algorithm: 

\newpage

\noindent
$G = - \frac{1}{f} A A^* Q$ \\
$do$ $i=r$, $1$, $-1$ \\
\hspace*{10mm} $t_1=b_i$ \\
\hspace*{10mm} $do$ $j=1$, $i-1$ \\
\hspace*{20mm} $z_j = q_j^* t_j$ \\
\hspace*{20mm} $t_{j+1} = t_j - z_j q_j$ \\
\hspace*{10mm} $enddo$ \\
\hspace*{10mm} $\displaystyle g_i=\frac{g_i-q_i q_i^* g_i}{||t_i||_2}$ \\
\hspace*{10mm} $do$ $j=i-1$, $1$, $-1$ \\
\hspace*{20mm} $\alpha = q_j^* g_i$ \\
\hspace*{20mm} $g_j = g_j - z_j^* g_i - \alpha^* t_j$ \\
\hspace*{20mm} $g_i = g_i - \alpha q_j$ \\
\hspace*{10mm} $enddo$ \\
$enddo$ \\

\noindent Here we use the $Z = (z_1, \dots, z_r) \in \C^r$ array with only
$r$ elements for better performance.

The total arithmetical complexity of the computation of the gradient is
$4mr(2r+n)$ operations. If we compare this with MGS ($2mr(r+n)$), it is
less than 4 times greater.

We obtain similar results for the Gramm-Schmidt (not MGS) orthogonalization:
it needs $\displaystyle 3mr + \frac{r(r+1)}{2}$ words in memory and
works with $2mr(3r+2n)$ operations, but because of stability issues we do not
recommend using it.

\bigskip

{\bf\large

\noindent
NUMERICAL EXPERIMENTS

}
\bigskip

First we compare the general characteristics of our new approach with
those of well-know approaches. We create the complex matrices $A$ and
$B$ with random numbers, compute derivatives for different
sizes of the problems by our new methods based on Gramm-Schmidt (AGS)
and Modified Gramm-Schmidt (AMGS) algorithms, and compare our methods with the finite
difference (FD) and Baur-Strassen (BS) methods (Tables 1, 2).

Furthermore, we show how those algorithms work. We perform a set of experiments
and check the number of iterations for convergence of the Broyden method
\cite{DShnab}. In this set of experiments, the matrix $B$ is a real matrix with
$b_i = p_i \otimes q_i \in \R^{n^2}$, $i=1, \dots, R$, where $\otimes$
is the Kronecker product of vectors and $p_i, q_i \in \R^n $ are unknown
vectors. We change $n \in [2, 20]$ and $r \in [2, 20]$ (Table 3).
This problem occurs in the three-way decomposition \cite{ENUMATH,NLAA}.

Hence, our new method (AMGS) is stable enough (like the BS method), yet also up to a thousand times 
faster than BS and FD methods and does not require much
additional memory (only 4 times more than FD).

\noindent
\underline{Table 1.} Memory requirements (in words) for FD, BS, AGS, and AMGS
methods.
$$\begin{tabular}{|c|c|c|c|c|c|c|} \hline
$m$ & $n$ & $r$ & FD    & BS    & AGS   & AMGS  \\ \hline
 2  &  1  &  1  & 4     &  152  & 12    &  14   \\ \hline
10  &  2  &  2  & 40    & 2.5k  & 92    & 124   \\ \hline
100 &  10 &  10 & 1.9k  & 619k  & 4.1k  & 5.9k  \\ \hline
1000& 100 & 100 & 195k  & 610m  & 410k  & 586k  \\ \hline \hline
1000& 100 & 10  & 19.5k & 33.5m & 39.3k & 58.6k \\ \hline
1000& 10  & 100 & 195k  & 337m  & 410k  & 586k  \\ \hline
\end{tabular}$$
\underline{Table 2.} Computational time of FD, BS, AGS, and AMGS methods.
$$\begin{tabular}{|c|c|c|c|c|c|c|} \hline
$m$ & $n$ & $r$ & FD        & BS        & AGS       & AMGS      \\ \hline
10  &  2  &  2  & $30$us    & $10$us    & $30$us    & $10$us    \\ \hline
100 &  10 &  10 & $209$ms   & $5.5$us   & $160$us   & $110$us   \\ \hline
1000& 100 & 100 & $9.5$h    & $6.4$s    & $175$ms   & $208$us   \\ \hline \hline
1000& 100 & 10  & $.92$h    & $359$us   & $9.2$ms   & $9.3$ms   \\ \hline
1000& 10  & 100 & $5$h      & $3.9$s    & $154$ms   & $205$ms   \\ \hline
\end{tabular}$$
\underline{Table 3.} The dependence of the total number of iterations in the
Broyden method on the method of gradient computation and problem size
for the first series of experiments ($N_u$ is the total number of unknowns).
$$\begin{tabular}{|c|c|c|c|c|c|c|c|} \hline
           $m$ & $n$  & $r$  & $N_u$    & FD      & BS     & AGS     & AMGS   \\ \hline
  $2 \times 2$ & $2$  & $2$  & $8$      & $4$     & $4$    & $4$     & $4$    \\ \hline
  $5 \times 5$ & $5$  & $5$  & $125$    & $244$   & $238$  & $521$   & $238$  \\ \hline
$10 \times 10$ & $10$ & $10$ & $1000$   & $3061$  & $1581$ & $>5000$ & $1581$ \\ \hline
$20 \times 20$ & $20$ & $20$ & $8000$   & $>5000$ & $2464$ & $>5000$ & $2464$ \\ \hline
\end{tabular}$$

\end{document}